\documentclass[10pt]{article}
\usepackage{amssymb}
\usepackage{amsmath}
\usepackage{verbatim}
\usepackage{pifont}
\usepackage{geometry}

\newtheorem{theorem}{Theorem}

\newtheorem{definition}[theorem]{Definition}

\input{tcilatex}

\begin{document}

\title{A majorization method for localizing graph topological indices }
\author{\emph{Monica Bianchi $^{1}$} \thanks{%
e-mail: monica.bianchi@unicatt.it} \and \emph{Alessandra Cornaro $^{1}$}
\thanks{%
Corresponding author, e-mail: alessandra.cornaro@unicatt.it,
Tel:+390272343683, Fax:+390272342324} \and \emph{Anna Torriero$^{1}$}
\thanks{%
e-mail: anna.torriero@unicatt.it} }
\date{ }
\maketitle

\begin{abstract}
This paper presents a unified approach for localizing some relevant graph
topological indices via majorization techniques. Through this method, old
and new bounds are derived and numerical examples are provided, showing how
former results in the literature could be improved.


\end{abstract}

\centerline{$^{{}^1}$ Department of Mathematics and Econometrics, Catholic University, Milan,
Italy.}

\section{Introduction and preliminaries}

\noindent Optimization problems involving Schur-convex or Schur-concave
functions have received much attention in the literature providing some
useful applications in different fields: for example they became a
convenient tool to localize eigenvalues of a real spectrum matrix (\cite%
{Bianchi1}, \cite{Tarazaga2}) or more generally to obtain bounds for an
arbitrary order statistic distribution (\cite{Boyd}). More recently some
issues related to the structural properties of graphs, characterized in
terms of invariants, have been explored solving suitable optimization
problems via majorization techniques (\cite{BCT}, \cite{Grassi}).

\noindent The aim of this paper is to explore this approach, showing that it
is a powerful tool for achieving, through a unified scheme, many well-known
bounds of some graph topological indices that can be expressed in terms of
Schur-convex or Schur-concave functions. One major advantage of this
technique is to set up a strategy for improving these bounds as well as to
obtain new ones. We start recalling some basic definitions on majorization.
More details can be found in \cite{Marshall}.


\begin{definition}
Given two vectors $\mathbf{y}$, $\mathbf{z\in }$ $D=\{\mathbf{x}\in \mathbb{R%
}^{n}:x_{1}\geq x_{2}\geq ...\geq x_{n}\}$, the majorization order $\mathbf{y%
}\trianglelefteq \mathbf{z}$ means:
\begin{equation*}
\begin{cases}
\left\langle \mathbf{y},\mathbf{s}^{\mathbf{k}}\right\rangle \leq
\left\langle \mathbf{z},\mathbf{s}^{\mathbf{k}}\right\rangle ,\text{ }%
k=1,...,(n-1) \\
\left\langle \mathbf{y},\mathbf{s}^{\mathbf{n}}\right\rangle =\left\langle
\mathbf{z},\mathbf{s}^{\mathbf{n}}\right\rangle%
\end{cases}%
\end{equation*}%
where $\left\langle \cdot ,\cdot \right\rangle $ is the inner product in $%
\mathbb{R}^{n}$ and $\mathbf{s^{j}}=[\underbrace{1,1,\cdots ,1}_{j},%
\underbrace{0,0,\cdots 0}_{n-j}],\quad j=1,2,\cdots ,n.$
\end{definition}

%

\begin{definition}
Given a set $S\subseteq D\cap \{\mathbf{x}\in \mathbb{R}^{n}:\left\langle
\mathbf{x},\mathbf{s^{n}}\right\rangle =a\}$, a vector $\mathbf{x}^{\ast
}(S)\in S$ ($\mathbf{x}_{\ast }(S)\in S$) is said to be the maximal
(minimal) vector in $S$ with respect to the majorization order if $\mathbf{x}%
\trianglelefteq \mathbf{x}^{\ast }(S)$ ($\mathbf{x}_{\ast
}(S)\trianglelefteq \mathbf{x}$) for each $\mathbf{x}\in S$.
\end{definition}

\noindent Since the upper and lower level sets:%
\begin{equation*}
U(\mathbf{z})=\left\{ \mathbf{x\in }\text{ }S_{a}:\mathbf{z\trianglelefteq x}%
\right\} ,\text{ }L(\mathbf{z})=\left\{ \mathbf{x\in }\text{ }S_{a}:\mathbf{%
x\trianglelefteq z}\right\}
\end{equation*}%
are closed, the existence of maximal and minimal elements of $S$ are ensured
by its compactness. In \cite{Bianchi1}, \cite{BCT} and \cite{Marshall}, the
extremal vectors for some closed subsets of
\begin{equation*}
\Sigma _{a}=\{\mathbf{x}\in \mathbb{R}^{n}:x_{1}\geq x_{2}\geq ...\geq
x_{n}\geq 0\}\cap \{\mathbf{x}\in \mathbb{R}^{n}:\left\langle \mathbf{x},%
\mathbf{s^{n}}\right\rangle =a\}
\end{equation*}%
have been computed. In particular in \cite{BCT} the maximal and minimal
elements of
\begin{equation*}
S_{a}=\Sigma _{a}\cap \{\mathbf{x}\in \mathbb{R}^{n}:M_{i}\geq x_{i}\geq
m_{i},i=1,\cdots ,n\}
\end{equation*}%
where $M_{1}\geq M_{2}\geq \cdots \geq M_{n},$ $m_{1}\geq m_{2},\cdots \geq
m_{n}$, have been derived, studying also the extremal vectors of some
subsets of $S_{a}$ of particular interest: 

\begin{equation}
S_{a}^{1}=\Sigma _{a}\cap \{\mathbf{x}\in \mathbb{R}^{n}:M\geq x_{1}\geq
x_{2}\geq \cdots \geq x_{n}\geq m\}\text{ ( see also \cite{Marshall})}
\label{S1}
\end{equation}%
\begin{equation}
S_{a}^{2}=\Sigma _{a}\cap \{\mathbf{x}\in \mathbb{R}^{n}:x_{i}\geq \alpha ,%
\text{ }i=1,\cdots ,h,\text{ }1\leq h\leq n,\text{ }\alpha \leq \frac{a}{h}\}%
\text{ ( see also \cite{Bianchi1})}  \label{S2}
\end{equation}%
\begin{equation}
S_{a}^{3}=\Sigma _{a}\cap \{\mathbf{x}\in \mathbb{R}^{n}:x_{i}\leq \alpha
,i=h+1,\cdots ,n,\text{ }1\leq h\leq (n-1),\text{ }\alpha <a\};\text{ ( see
also \cite{Bianchi1})}  \label{S3}
\end{equation}%
\begin{eqnarray}
S_{a}^{(h)} &=&\Sigma _{a}\cap \left\{ \mathbf{x}\in \mathbb{R}%
^{n}:M_{1}\geq x_{1}\geq \cdots \geq x_{h}\geq m_{1},\right.  \label{S4} \\
&&\left. M_{2}\geq x_{h+1}\geq \cdots \geq x_{n}\geq m_{2},\text{ }%
m_{i}<M_{i},i=1,2\right\} .  \notag
\end{eqnarray}

\begin{definition}
A symmetric function $\phi $: $A\rightarrow
\mathbb{R}
$, $A\subseteq
\mathbb{R}
^{n}$, is said to be Schur-convex on $A$ if $\mathbf{x}\trianglelefteq
\mathbf{y}$ implies $\phi (\mathbf{x})\leq $ $\phi (\mathbf{y})$. If in
addition $\phi (\mathbf{x})<$$\phi (\mathbf{y})$ for $\mathbf{x}%
\trianglelefteq \mathbf{y}$ but $\mathbf{x}$ is not a permutation of $%
\mathbf{y}$, $\phi $ is said to be strictly Schur-convex on $A$. A function $%
\phi $ is (strictly) Schur-concave on $A$ if $-\phi $ is (strictly)
Schur-convex on $A.$
\end{definition}

\noindent Given an interval $I\subset
\mathbb{R}
$, and a (strictly) convex function $g:I\rightarrow
\mathbb{R}
$, the function $\phi (\mathbf{x})=\sum_{i=1}^{n}g(x_{i})$ is (strictly)
Schur-convex on $I^{n}=\underbrace{I\times I \times \cdots \times I}_{n-
times}$. The corresponding result holds if $g$ is (strictly) concave on $I^n$%
.

\noindent From the order preserving property of Schur-convex functions, the
solution of some constrained nonlinear optimization problems of particular
interest can be obtained in a straightforward way.
More precisely, the problem we face is the following:
\begin{equation}
\begin{cases}
\max \,\,\,(\min )\text{ }\phi (\mathbf{x})\tag{P}\label{P} \\
\text{subject to }\mathbf{x}\in S%
\end{cases}%
\end{equation}%
\noindent where $S\subseteq \mathbb{R}^{n}$ is a generic subset which admits
maximal vector $\mathbf{x}^{\ast }(S)$ and minimal vector\textit{\ }$\mathbf{%
x}_{\ast }(S)$ with respect to the majorization order. If the objective
function $\phi $\textit{\ }is Schur-convex, the maximum and the minimum are
attained in $\mathbf{x}^{\ast }(S)$\textbf{\textit{\ }}and $\mathbf{x}_{\ast
}(S)$ respectively; the opposite holds if $\phi $ is a Schur-concave
function. This enable us to solve problem (P) in a more direct way, avoiding
the standard approach of Lagrange multipliers. In the next section we
illustrate how the procedure above can be successfully applied to get upper
and lower bounds for some relevant topological indices, which found
applications in many fields varying from chemistry to network analysis.


\section{Majorization and bounds for graph topological indices.}

\noindent Let us firstly recall some basic graph notations and concepts to
be used later. For more details refer to \cite{Grone_Merris} and \cite%
{Harary}. \bigskip Let $G=(V,E)$ be a simple, connected, undirected graph
where $V=\left\{ v_{1},...,v_{n}\right\} $ is the set of vertices and $%
E\subseteq V\times V$ the set of edges. We consider graphs with fixed order $%
\left\vert V\right\vert =n\ $ and fixed size $\left\vert E\right\vert =m.$
Denote by $\pi =(d_{1},d_{2},..,d_{n})$ the degree sequence of $G,$ where $%
d_{i}$ is the degree of vertex $v_{i}$, arranged in non increasing order $%
d_{1}\geq d_{2}\geq \cdots \geq d_{n}$. It is well known that $\overset{n}{%
\underset{i=1}{\dsum }}d_{i}=2m$ and that if $G$ is a tree, i.e. a connected
graph without cycles, $m=n-1.$ Let $A(G)$ be the adjacency matrix of $G$ and
$\lambda _{1}\geq \lambda _{2}\geq ...\geq \lambda _{n}$ be the set of
(real) eigenvalues of $A(G).$ The matrix $L(G)=D(G)-A(G)$ is called
Laplacian matrix of $G$, where $D(G)$ is the diagonal matrix of vertex
degrees. Let $\mu _{1}\geq \mu _{2}\geq ...\geq \mu _{n}$ be the eigenvalues
of $L(G)$. The following properties of spectra of $A(G)$ and $L(G)$ are well
known :

\begin{equation*}
\overset{n}{\underset{i=1}{\dsum }}\lambda _{i}=0;\,\,\,\overset{n}{\underset%
{i=1}{\dsum }}\lambda _{i}^{2}=2m;
\end{equation*}%
\begin{equation*}
\overset{n}{\underset{i=1}{\dsum }}\mu _{i}=2m;\,\,\,\mu _{1}\geq
1+d_{1}\geq \dfrac{2m}{n};\,\,\,\mu _{n}=0,\text{ }\mu _{n-1}>0.
\end{equation*}

\noindent Notice that most of the topological indices of graphs are
formulated by strictly Schur-convex (Schur-concave) functions of the degree
sequence as well as the eigenvalues of $A(G)$ or $L(G)$. The corresponding
bounds are generally expressed in terms of size and order of $G$ but they
can take also into account the degrees of one or more vertices of $G.$ With
respect to the degree sequence, one of the most popular index is the General
Randic index:
\begin{equation*}
R_{\alpha }(G)\underset{\left( v_{i},v_{j}\right) \in E}{=\sum \left(
d_{i}d_{j}\right) ^{\alpha }},
\end{equation*}%
where $\alpha $ is a non zero real number (see \cite{Bollobas}). Setting a
specific value for $\alpha $, some very well known indices can be obtained:
for example, $\alpha =1$ corresponds to the Zagreb index $M_{2}(G)$ (see
\cite{Nikolic}) while $\alpha =-\dfrac{1}{2}$ and $\alpha =-1$ to the
branching indices (see \cite{Randic}). In the recent paper \cite{Zhou-Trina}%
, the sum-connectivity index has been proposed and in \cite{DuZhouTrina} it
has been extended to the Generalized sum-connectivity index defined as
follows:
\begin{equation*}
\chi _{\alpha }(G)=\sum_{(v_{i},v_{j})\in E}\left( d_{i}+d_{j}\right)
^{\alpha }.
\end{equation*}%
Note that for $\alpha =1$ we obtain the first Zagreb index $%
M_{1}(G)=\sum_{(v_{i},v_{j})\in E}\left( d_{i}+d_{j}\right) =\sum_{v_{i}\in
E}d_{i}^{2}.$ Other frequently used indices involve Schur-convex or
Schur-concave functions of the eigenvalues of $A(G)$ and $L(G)$. We recall,
among the others:

\begin{enumerate}
\item Energy index: $E(G)=\sum_{i=1}^{n}|\lambda _{i}|$ (\cite{Gutman_Energy}%
)

\item $s_{\alpha }(G)=\sum_{i=1}^{n-1}\mu _{i}^{\alpha },\alpha \neq 0,1$ (%
\cite{LiuLiu}, \cite{Zhou})

\item Kirchhoff index $Kf(G)=n\sum_{i=1}^{n-1}\dfrac{1}{\mu _{i}}=ns_{-1}(G)$
(\cite{ZhouTrina1}, \cite{ZhouTrina2})

\item Laplacian Estrada index $LEE(G)=\sum_{i=1}^{n}e^{\mu _{i}}$ (\cite%
{ZhouGutman}, \cite{Zhou_estrada})
\end{enumerate}

\noindent In the next section we present theoretical results on upper and
lower bounds for the indices previously mentioned, recovering well known
bounds and improving, in some cases, the existing ones of the literature. We
postpone to Section3 some numerical examples.


\subsection{General Randic index}

\noindent The generalized Randic index can be equivalently expressed as:
\begin{equation*}
R_{\alpha }(G)=\underset{\left( v_{i},v_{j}\right) \in E}{\sum }\left(
d_{i}d_{j}\right) ^{\alpha }=\frac{1}{2}\left( \underset{\left(
v_{i},v_{j}\right) \in E}{\sum }\left( d_{i}^{\alpha }+d_{j}^{\alpha
}\right) ^{2}-\overset{n}{\underset{i=1}{\dsum }}d_{i}^{2\alpha +1}\right) .
\end{equation*}

\noindent Let $\pi =(d_{1},d_{2},..,d_{n})$ be a fixed degree sequence and $%
\mathbf{x}\in
\mathbb{R}
^{m}$ be the vector whose components are $d_{i}^{\alpha }+d_{j}^{\alpha },$
with $\left( v_{i},v_{j}\right) \in E.$ Hence, taking into account that $%
\overset{n}{\underset{i=1}{\sum }}d_{i}^{2\alpha +1}$ is a constant, $%
R_{\alpha }(G)$ is a Schur convex function of the degree sequence of $G$ and
it is minimal (maximal) if and only $f(\mathbf{x})=\overset{m}{\underset{i=1}%
{\sum }}x_{i}^{2}=$ $\left\Vert x\right\Vert _{2}^{2}$ is minimal(maximal).
Extending the result proved in \cite{Lovasz}, it is possible to show that
\begin{equation}
\sum_{i=1}^{m}x_{i}=\underset{\left( v_{i},v_{j}\right) \in E}{\sum }\left(
d_{i}^{\alpha }+d_{j}^{\alpha }\right) =\sum_{i=1}^{n}d_{i}^{\alpha +1}
\end{equation}%
and thus $\sum_{i=1}^{m}x_{i}$ is a constant. Hence, considering a subset $S$
of $\Sigma _{a}\subseteq \mathbb{R}^{m},$ where $a=\sum_{i=1}^{n}d_{i}^{%
\alpha +1}$, 
which admits $x_{\ast }(S)$ and $x^{\ast }(S)$ as extremal vectors with
respect to the majorization order, the function $f$ attains its minimum and
maximum on $S$ at $f(x_{\ast }(S))$ and $f(x^{\ast }(S))$, respectively. The
general Randic index can be consequently bounded as follows:

\begin{equation}
\frac{\left\Vert x_{\ast }(S)\right\Vert _{2}^{2}-\overset{n}{\underset{i=1}{%
\sum }}d_{i}^{2\alpha +1}}{2}\leq R_{\alpha }(G)\leq \frac{\left\Vert
x^{\ast }(S)\right\Vert _{2}^{2}-\overset{n}{\underset{i=1}{\sum }}%
d_{i}^{2\alpha +1}}{2}.  \label{randic}
\end{equation}%
Clearly different numerical bounds can be derived, if we characterize
suitably the set $S$. The structure of $S$, among the typology of sets given
in (\ref{S1})-(\ref{S4}), will depend on the information available on the
degree sequence of $G.$

\noindent For the particular case $\alpha =1$, which corresponds to the
Zagreb index $M_{2}$, this methodology was applied in \cite{Grassi},
and, more recently, in \cite{BCT} where the authors get sharper bounds for
the index $M_{2}$, valid for a particular class of graphs having exactly $h$
pendant vertices, i.e. vertices with degree one, with degree sequence
\begin{equation}
\pi =(d_{1},\cdots ,d_{n-h},\underbrace{1,\cdots ,1}_{h}).  \label{ds}
\end{equation}

\noindent Starting from (6), the following results allow to compute the
Randic index corresponding to a generic $\alpha ,$ for a class of graphs
with $h$ pendant vertices and degree sequence given by (\ref{ds}). Let us
assume, furthermore, that the degree sequence $\pi $ is such that:
\begin{equation}
\begin{cases}
1+d_{1}^{\alpha }\leq d_{n-h}^{\alpha }+d_{n-h-1}^{\alpha } & \text{ for }%
\alpha >0 \\
d_{n-h}^{\alpha }+d_{n-h-1}^{\alpha }\leq 1+d_{1}^{\alpha } & \text{ for }%
\alpha <0%
\end{cases}%
.  \label{degree_ass}
\end{equation}

\noindent These conditions, always satisfied for $\alpha =-1$, assure that
the $h$ components of $\mathbf{x}$\textbf{,} corresponding to pendant nodes,
are separated from the others. Moreover, since the components of $\mathbf{x}$
are arranged in nonincreasing order, for $\alpha >0,$ they are assigned to
the last $h$ positions of the vector $\mathbf{x}$ while, for $\alpha <0$, to
the first $h$ positions of the vector $\mathbf{x}$. Being interested in
bounds for the branching indices associated to $\alpha =-\frac{1}{2}$ and $%
\alpha =-1$, in the sequel we only focus on negative values of $\alpha $. It
is easy to see that the choice of a degree sequence equal to (\ref{ds}),
combined with (\ref{degree_ass}), lead us to face the following set of type (%
\ref{S4}):

\begin{equation}
\begin{split}
S=\{\mathbf{x}\in \mathbb{R}^{m}:& \,\,\,1+d_{1}^{\alpha }\leq x_{h}\leq
\cdots \leq x_{1}\leq d_{n-h}^{\alpha }+1, \\
& d_{1}^{\alpha }+d_{2}^{\alpha }\leq x_{m}\leq \cdots \leq x_{h+1}\leq
d_{n-h}^{\alpha }+d_{n-h-1}^{\alpha },\text{ }\left\langle \mathbf{x},%
\mathbf{s}^{m}\right\rangle =\sum_{i=1}^{n}d_{i}^{\alpha +1}\}
\end{split}
\label{S4n}
\end{equation}%
whose maximal and minimal elements with respect to the majorization order
can be computed by corollaries 3 and 10 in \cite{BCT}. Let us notice that
both inequalities in (\ref{randic}) are attained if and only if the set $S$
reduces to a singleton. This is the case when $d_{1}=d_{n-h}$, i.e. when all
non pendant vertices have the same degree (see \cite{BCT} for some
significant examples).

\subsection{Generalized sum-connectivity index}

\noindent Let $\pi $ be a fixed degree sequence and $\mathbf{x}\in \mathbb{R}%
^{m}$ be the vector whose components are $(d_{i}+d_{j})$, $(v_{i},v_{j})\in
E.$ The function $f(\mathbf{x})=\sum_{i=1}^{m}x_{i}^{\alpha }$ is strictly
Schur-convex for $\alpha >1\text{ or }\alpha <0$, while it is strictly
Schur-concave for $0<\alpha <1$. Thus, taking into account that $%
\sum_{i=1}^{m}x_{i}=\sum_{i=1}^{n}d_{i}^{2}$ is a constant, considering a
subset $S$ of $\Sigma _{a}$, where $a=\sum_{i=1}^{n}d_{i}^{2}$, for $\alpha
>1\text{ or }\alpha <0$ we get
\begin{equation}
\left\Vert x_{\ast }(S)\right\Vert _{\alpha }^{\alpha }\leq \chi _{\alpha
}(G)\leq \left\Vert x^{\ast }(S)\right\Vert _{\alpha }^{\alpha },
\label{sum_con}
\end{equation}%
where $\left\Vert \cdot \right\Vert _{\alpha }$ stands for the $l_{\alpha }-$%
norm. For $0<\alpha <1$, the bounds are exchanged. 



\subsection{Energy index}


\noindent Let us point out that $E(G)=\sum_{i=1}^{n}|\lambda
_{i}|=\sum_{i=1}^{n}\sqrt{\lambda _{i}^{2}}$ and thus this index is a
Schur-concave function of the variables $\lambda _{i}^{2},i=1,\cdots ,n.$
Furthermore $\sum_{i=1}^{n}\lambda _{i}^{2}=2m$. Let $x_{i}=\lambda
_{i}^{2},i=1,\cdots ,n.$ It is well known that $\lambda _{1}\geq \frac{2m}{n}
$. If a sharper lower bound for $\lambda _{1}$ is available, i.e. $\lambda
_{1}\geq k(\geq \frac{2m}{n})$, we have $x_{1}\geq k^{2}\geq \left( \frac{2m%
}{n}\right) ^{2}\geq \frac{2m}{n}$. The minimal element of the set
\begin{equation*}
S=\Sigma _{a}\cap \{\mathbf{x}\in \mathbb{R}^{n}:x_{1}\geq k^{2}\}
\end{equation*}%
where $a=2m,$ is given, using Corollary 14 in \cite{BCT}, by:
\begin{equation*}
x_{\ast }(S)=\left[ k^{2},\underbrace{\dfrac{2m-k^{2}}{n-1},\cdots ,\dfrac{%
2m-k^{2}}{n-1}}_{n-1}\right] .
\end{equation*}

\noindent Thus
\begin{equation*}
E(G)\leq k+\sqrt{(n-1)\left( 2m-k^{2}\right) }.
\end{equation*}%
In the case of bipartite graphs, since $\lambda _{1}=-\lambda _{n}$, we have
$x_{1}=x_{2}$ and thus we face the set
\begin{equation*}
S=\Sigma _{a}\cap \{\mathbf{x}\in \mathbb{R}^{n}:x_{i}\geq k^{2},i=1,2\}
\end{equation*}%
whose minimal element becomes
\begin{equation*}
x_{\ast }(S)=\left[ k^{2},k^{2},\underbrace{\dfrac{2m-2k^{2}}{n-2},\cdots ,%
\dfrac{2m-2k^{2}}{n-2}}_{n-2}\right]
\end{equation*}%
and as a consequence
\begin{equation*}
E(G)\leq 2k+\sqrt{(n-2)\left( 2m-2k^{2}\right) }.
\end{equation*}

\noindent For $k=\frac{2m}{n}$ we get the bounds in \cite{Koolen1} and \cite%
{Koolen2}, while for $k=\frac{\sum_{i=1}^{n}d_{i}^{2}}{n}$ the bounds in
\cite{Zhou_energy}. Different choices of $k$ allow to obtain all the other
bounds in \cite{Liu_Lu_Tan_energy} and a new bound for $E(G)$ as soon as a
sharper lower bound of $\lambda _{1}$ is available.


\subsection{Laplacian indices}


\begin{enumerate}
\item Let $s_{\alpha }(G)=\sum_{i=1}^{n-1}\mu _{i}^{\alpha },\alpha \neq 0,1$
be the index given by the sum of the $\alpha $-th power of the non zero
Laplacian eigenvalues and consider the set
\begin{equation*}
S=\Sigma _{a}\cap \{\mathbf{\mu }\in \mathbb{R}^{n-1}:\mu _{1}\geq 1+d_{1}\}
\end{equation*}%
where $a=\sum_{i=1}^{n-1}\mu _{i}=2m$ and $\Sigma _{a}\subseteq \mathbb{R}%
^{n-1}$. Since $\frac{2m}{n-1}\leq (1+d_{1})$, by Corollary 14 in \cite{BCT}
it follows that the minimal element of $S$ is
\begin{equation*}
x_{\ast }(S)=\left[ 1+d_{1},\underbrace{\frac{2m-1-d_{1}}{n-2},\cdots ,\frac{%
2m-1-d_{1}}{n-2}}_{n-2}\right] .
\end{equation*}%
Taking into account the Schur-convexity or Schur-concavity of the functions $%
s_{\alpha }(G)$, the bounds in \cite{Zhou}, Theorem 3 can be easily derived.
The same approach can be used to find the bounds for bipartite graphs given
in \cite{Zhou}, Theorem 5, observing that in this case $\mu _{1}\geq 2\sqrt{%
\frac{\sum_{i=1}^{n}d_{i}^{2}}{n}}$.

The bounds in \cite{Zhou}, Theorem 3 can be also improved taking into
account further information on the localization of the eigenvalues. For
instance, since $\mu _{2}\geq d_{2}$ (see \cite{BroHae}), we can consider
the set
\begin{equation}
S=\Sigma _{a}\cap \{\mathbf{\mu }\in \mathbb{R}^{n-1}:\mu _{1}\geq
1+d_{1},\,\,\mu _{2}\geq d_{2}\}  \label{secondolaplaciano}
\end{equation}%
whose minimal element, computed applying Theorem 8 in \cite{BCT}, is given
by
\begin{equation*}
x_{\ast }(S)=\left[ 1+d_{1},d_{2},\underbrace{\frac{2m-1-d_{1}-d_{2}}{n-3}%
,\cdots ,\frac{2m-1-d_{1}-d_{2}}{n-3}}_{n-3}\right]
\end{equation*}%
Thus
\begin{equation}
\begin{cases}
s_{\alpha }(G)\geq (1+d_{1})^{\alpha }+d_{2}^{\alpha }+\frac{%
(2m-1-d_{1}-d_{2})^{\alpha }}{(n-3)^{\alpha -1}} & \alpha >1,\alpha <0 \\
s_{\alpha }(G)\leq (1+d_{1})^{\alpha }+d_{2}^{\alpha }+\frac{%
(2m-1-d_{1}-d_{2})^{\alpha }}{(n-3)^{\alpha -1}} & 0<\alpha <1%
\end{cases}%
.
\end{equation}


\item Notice that bounds on the Kirchhoff index $Kf(G)=n\sum_{i=1}^{n-1}%
\dfrac{1}{\mu _{i}}=ns_{-1}(G)$ can be easily derived by bounds on $%
s_{-1}(G) $. In particular, using the set (\ref{secondolaplaciano}) we get
\begin{equation}
Kf(G)\geq n\left( \frac{1}{1+d_{1}}+\frac{1}{d_{2}}+\frac{(n-3)^{2}}{%
(2m-1-d_{1}-d_{2})}\right) .  \label{KG}
\end{equation}

\item[3.] The Laplacian Estrada index $LEE(G)=\sum_{i=1}^{n}e^{\mu _{i}}$ is
a Schur-convex function of the Laplacian eigenvalues. Thus, proceeding as
before we can easily recover the bounds in \cite{Zhou_estrada}, Remark 2.
Furthermore, if we use the set (\ref{secondolaplaciano}), we find the bound:%
\begin{equation}
LEE(G)\geq 1+e^{1+d_{1}}+e^{d_{2}}+(n-3)e^{\frac{2m-1-d_{1}-d_{2}}{n-3}}.
\label{LEE}
\end{equation}
\end{enumerate}

\section{Numerical examples}

In this section we provide some numerical examples to illustrate our results
and to show that, at least for the considered cases, we achieve sharper
bounds with respect to recent literature.

\noindent

\begin{enumerate}
\item[i)] \textbf{General Randic index}
\end{enumerate}

\noindent Let us consider a tree $T$ with the degree sequence $\pi
=(5,3,3,3,3,3,1,1,1,1,1,1,1,1,1,1)$. We start exploring the Randic index $%
R(T)$, with $\alpha =-\dfrac{1}{2}.$ The comparison (see Table 1) with
bounds in \cite{Li_Shi} shows that bounds obtained by applying Corollary 3
and 10 in \cite{BCT} always perform better:

\begin{table}[th]
\centering
\begin{tabular}{|l|l|l|}
\hline
Ref. & Lower & Upper \\ \hline
\cite{Li_Shi} & $5.891\,249$ (see Theorem 2.6) & $7.486\,884$ (see Theorem
2.24) \\ \hline
\cite{BCT} & $6.768\,317$ & $7.064\,497$ \\ \hline
\end{tabular}%
\caption{Lower and upper bounds for $R(T)$. }
\end{table}
\newpage \noindent For the case with $\alpha =-1$ and the tree $T$ with the
same degree sequence, making a comparison with the result in \cite{Li_Shi},
Theorem 3.7, we get

\begin{table}[th]
\centering%
\begin{tabular}{|l|l|l|}
\hline
Ref & Lower & Upper \\ \hline
\cite{Li_Shi} & $1$ & $4.888889$ \\ \hline
\cite{BCT} & $3.2$ & $3.666667$ \\ \hline
\end{tabular}%
\caption{Lower and upper bounds for $R_{-1}(T)$. }
\end{table}

\noindent

\begin{enumerate}
\item[ii)] \textbf{Generalized sum-connectivity index}
\end{enumerate}

\noindent Let us consider a graph with degree sequence $\pi =\left(
3,2,2,1\right) $ and pick $\alpha =-5.$ Since, in this case, $n$ is even, we
make a comparison (see Table 3) between the results in \cite{DuZhouTrina}
and those obtained from Corollary 3 in \cite{BCT}.
\begin{table}[th]
\centering
\begin{tabular}{|l|l|}
\hline
Ref. & Upper \\ \hline
\cite{DuZhouTrina} & $9.207015\times 10^{-3}$ (see Theorem 1) \\ \hline
\cite{BCT} & $5.075226\times 10^{-3}$ \\ \hline
\end{tabular}%
\caption{Upper bounds for $\protect\chi _{\protect\alpha }(G)$ with $n$
even. }
\end{table}
\noindent We can present here also the case when $n$ is odd, always choosing
$\alpha =-5$. \noindent For a tree with the following degree sequence $\pi
=\left( 4,4,4,3,1,1,1,1,1,1,1,1,1\right) $, computing the bounds (see Table
4) we have:
\begin{table}[!h]
\centering
\begin{tabular}{|l|l|}
\hline
Ref. & Upper \\ \hline
\cite{DuZhouTrina} & $2.487446\times 10^{-2}$ (see Theorem 1) \\ \hline
\cite{BCT} & $4.94124\times 10^{-3}$ (see Corollary 3) \\ \hline
\end{tabular}%
\caption{Upper bounds for $\protect\chi _{\protect\alpha }(G)$ with $n$ odd.}
\end{table}

\begin{enumerate}
\item[iii)] \textbf{Laplacian indices}
\end{enumerate}

\noindent Let us consider a graph whose degree sequence is $\pi
=(3,2,2,2,1). $

\begin{enumerate}
\item[a)] \textit{Kirchhoff index}
\end{enumerate}

For this index, making a comparison (see Table 5) between the result
obtained by applying the lower bound obtained in equation (\ref{KG}) and the
one derived from \cite{ZhouTrina2}, Proposition 1, we have
\begin{table}[!h]
\centering%
\begin{tabular}{|l|l|}
\hline
Ref & Lower \\ \hline
\cite{ZhouTrina2} & $12.5$ \\ \hline
formula (\ref{KG}) & $13.75$ \\ \hline
\end{tabular}%
\caption{Lower bounds for $Kf(G)$.}
\end{table}
\newpage

\begin{enumerate}
\item[b)] \textit{Laplacian Estrada index}
\end{enumerate}

In this case, we draw a comparison between the application of formula (\ref%
{LEE}) and the result provided in \cite{ZhouGutman}, Proposition 3.4 (see
Table 6)
\begin{table}[h]
\centering%
\begin{tabular}{|l|l|}
\hline
Ref & Lower \\ \hline
\cite{ZhouGutman} & $59.90761$ \\ \hline
formula (\ref{LEE}) & $68.42377$ \\ \hline
\end{tabular}%
\caption{Lower bounds for $LEE(G)$.}
\end{table}

\section{Conclusion}

In this paper we present a unified approach for localizing some relevant
graph topological indices, based on the optimization of \ Schur-convex or
Schur-concave functions. Our results have been derived through the
characterization of extremal vectors with respect to the majorization order,
under suitable constraints. We have shown that classical results can be
recovered and sometimes improved. Our theoretical approach paves the way to
find new bounds, by taking advantage of further information that can be
extracted from the graph $G$, that allows us to define tighter sets $S$ for
performing the optimization problem. Finally, other topological indices can
be localized whenever they can be expressed as Schur-convex or Schur-concave
functions.

\end{document}